\documentclass[11pt]{article}
\usepackage{}

\usepackage{amssymb}   
\usepackage{amsthm}    
\usepackage{amsmath}   
\usepackage{stmaryrd}  
\usepackage{titletoc}  
\usepackage{mathrsfs}  
\usepackage{graphicx}

\DeclareMathOperator*{\esssup}{ess\,sup}

\newlength{\defbaselineskip}
\newcommand{\setlinespacing}[1]%
           {\setlength{\baselineskip}{#1 \defbaselineskip}}

\theoremstyle{plain}
\newtheorem{thm}{Theorem}[section]
\newtheorem{cor}[thm]{Corollary}
\newtheorem{lem}[thm]{Lemma}
\newtheorem{prop}[thm]{Proposition}

\theoremstyle{definition}
\newtheorem{defn}{Definition}[section]

\newtheorem{rmk}{Remark}[section]

\newcommand{\eps}{\varepsilon}
\newcommand{\vf}{\varphi}

\newcommand{\bF}{\mathbb{F}}
\newcommand{\bS}{\mathbb{S}}

\newcommand{\bI}{\mathbb{I}}
\newcommand{\bR}{\mathbb{R}}

\newcommand{\bL}{\mathbb{L}}
\newcommand{\bP}{\mathbb{P}}
\newcommand{\sF}{\mathscr{F}}
\newcommand{\sH}{\mathscr{H}}
\newcommand{\sL}{\mathscr{L}}

\newcommand{\sB}{\mathscr{B}}
\newcommand{\sS}{\mathscr{S}}
\newcommand{\sU}{\mathscr{U}}
\newcommand{\la}{\langle}
\newcommand{\ra}{\rangle}

\newcommand{\rrow}{\rightarrow}

\makeatletter\@addtoreset{equation}{section} \makeatother
\begin{document}

\title{Neutral Backward Stochastic Functional Differential Equations and
Their Application}
\author{Wenning Wei\footnotemark[1] }

\footnotetext[1]{Department of Finance and Control Sciences, School
of Mathematical Sciences, and Laboratory of Mathematics for
Nonlinear Science, Fudan University, Shanghai 200433, China.
\textit{E-mail}: \texttt{071018033@fudan.edu.cn}.}

\maketitle

\begin{abstract}
In this paper we are concerned with a new type of backward equations
with anticipation which we call neutral backward stochastic
functional differential equations. We obtain the existence and
uniqueness and prove a comparison theorem. As an application, we
discuss the optimal control of neutral stochastic functional
differential equations, establish a Pontryagin maximum principle,
and give an explicit optimal value for the linear optimal control.

\end{abstract}

\textbf{Keywords}: Backward stochastic differential equations,
neutral stochastic functional differential equation, neutral
backward stochastic functional differential equations, optimal
control, duality, Pontryagin maximum principle

\textbf{Mathematical Subject Classification (2010)}: 60H10, 60H20,
93E20
\section{Introduction}
Throughout this paper, we fix $\delta>0$ as a positive constant. Let
$(\Omega,\sF,\bF,\bP)$ be a complete filtered probability space on
which a $d\textmd{-}$dimensional Brownian motion $\{W(t)\}_{t\geq0}$
is defined with $\{\sF_t\}_{t\geq0}$ being its natural filtration
augmented by all the $\bP\textmd{-}$null sets in $\sF$. Also we
define $\sF_u:=\sF_0$ for all $u\in[-\delta,0]$, then
$\bF:=\{\sF_t\}_{t\geq-\delta}$ is a filtration satisfying the usual
conditions on $[-\delta,+\infty)$.

In this paper, we investigate the following backward stochastic
equation with anticipation,
\begin{equation}\label{equ}
\left\{
\begin{split}
&-d\big[Y(t)-G(t,Y_t,Z_t)\big]=f(t,Y_t,Z_t)\,dt-Z(t)\,dW_t,~~t\in[0,T];\\
&Y(t)=\xi(t),~~~~~~t\in[T,T+\delta];\\
&Z(t)=\zeta(t),~~~~~~t\in[T,T+\delta],
\end{split}
\right.
\end{equation}
where $(Y_t,Z_t)$ denotes the path of the unknown processes $(Y,Z)$
on $[t,t+\delta]$, $G,f:
[0,T]\times\Omega\times\bL^2(0,\delta;\bR^n)\times\bL^2(0,\delta;\bR^{n\times
d})\rrow\bR^n$ are given maps, and $(\xi,\zeta)$ are given adapted
stochastic processes on $[T,T+\delta]$. We call $(G,f)$ the
generator. Equation \eqref{equ} is referred to as a neutral backward
stochastic functional differential equation (NBSFDE). It includes
many interesting cases.

When $G=0$ and $\delta=0$, \eqref{equ} becomes the well-known
backward stochastic differential equation (BSDE):
$$Y(t)=\xi(T)+\int_t^Tf(s,Y(s),Z(s))\,ds+\int_t^TZ(s)\,dW(s),~~t\in[0,T],$$
which was first introduced by Bismut \cite{Bismut73} for the linear
case and then extended to the nonlinear case by Pardoux and Peng
\cite{ParPeng_90}. It has been extensively applied in mathematical
finance and stochastic optimal control. For details, see El Karoui
et. al \cite{Karoui_Peng_Quenez97, El_Karoui01}, Peng \cite{Peng04},
Yong and Zhou \cite{YongBook99} and the references therein.

When $G=0$ and $f(t,\cdot,\cdot)$ only depends on the value of path
$(Y,Z)$ at t and $t+\delta$, Eq.\eqref{equ} becomes
\begin{equation*}
\left\{
\begin{split}
&-dY(t)=f(t,Y(t),Y(t+\delta),Z(t),Z(t+\delta))\,dt-Z(t)\,dW_t,~~t\in[0,T];\\
&Y(t)=\xi(t),~~~~~~t\in[T,T+\delta];\\
&Z(t)=\zeta(t),~~~~~~t\in[T,T+\delta].
\end{split}
\right.
\end{equation*}
It is the so-called anticipated backward stochastic differential
equations introduced by Peng and Yang \cite{PengYang09} when they
discuss optimal control of delayed stochastic functional
differential equations.

By an It\^{o} type neutral stochastic functional differential
equation (NSFDE), we mean the following:
\begin{equation*}
\left\{
\begin{split}
&d\big[X(t)-g(t,X^t)\big]=b(t,X^t)\,dt+\sigma(t,X^t)\,dW(t),~~t\in[0,T];\\
&X(t)=\vf(t),~~~t\in[-\delta,0],
\end{split}
\right.
\end{equation*}
where $X^t$ denotes the path of X on $[t-\delta,t]$. It is a type of
retarded functional equations. This equation was first introduced by
Kolmanovskii and Nosov \cite{Kolmanovskii_Nosov81} to model the
chemical engineering system. Since then, many papers are devoted to
the stability of the solutions. See Mao \cite{Mao95}, Huang and Mao
\cite{MaoX09},  and Randjelovi\'{c} and Jankovi\'{c}
\cite{Randjelovic07} and the references therein.

In deterministic case, optimal control of neutral functional
equations was discussed in 1960s and 1970s by Kolmanovskii and
Khvilon \cite{Kolmanovskii69}, Kent \cite{Kent71} and Banks and Kent
\cite{Banks72}. Optimal control of NSFDEs seems to remain to be
open. In this paper, via constructing the duality between linear
backward and forward stochastic neutral functional differential
equations, we discuss a simple optimal control in the stochastic
case. The more general one will be discussed elsewhere.

The rest of the paper is organized as follows. In section 2, we
prove the existence and uniqueness of adapted solutions and some
estimates of NBSFDEs. Section 3 is devoted to a comparison theorem.
Finally, as applications, we discuss an optimal control of NSFDEs,
construct a Pontryagin maximum principle, and obtain an explicit
optimal value of a linear optimal control.

\section{Existence and Uniqueness Result}
In this section, we prove the existence and uniqueness for NBSFDE
\eqref{equ}, construct some estimates of the solutions, and discuss
how $\delta$ affects the solutions in a simple case.

First, let us introduce some spaces. Let $H$ be a finite-dimensional
space like $\bR^n,\bR^{n\times d}, etc.$, whose norm is denoted by
$|\cdot|$. Denote by $\sB(D)$ the Borel $\sigma$-algebra of some
metric space D.
\begin{equation*}
\begin{split}
&\bL^2(s,\tau;H):=\left\{\psi:[s,\tau]\rrow
H~\Bigm|~\int_s^{\tau}|\psi(t)|^2\,dt<+\infty\right\},\\
&\bS^2(s,\tau;H):=\Big\{\psi:[s,\tau]\rrow H~|~\sup_{s\leq
t\leq\tau}|\psi(t)|^2<+\infty\Big\},\\
&\bL^2(\sF_T;H):=\Big\{\xi:\Omega\rrow
H~|~\sF_T\textmd{-}\textrm{measurable},~E[|\xi|^2]<+\infty\Big\},~~~~~~~~~~~~~~~~~~~~~~~~~~\\
\end{split}
\end{equation*}
\begin{equation*}
\begin{split}
&\sL^2_{\bF}(s,\tau;H):=\left\{X(\cdot):[s,\tau]\times\Omega\rrow
    H~\Bigm|~\bF\textmd{-}\textrm{adapted}~\textrm{and}~E\Big[\int_s^{\tau}|X(t)|^2\,dt\Big]<+\infty\right\},\\
&\sS^2_{\bF}([s,\tau];H):=\Big\{X(\cdot):[s,\tau]\times\Omega\rrow
H~|~\bF\textmd{-}\textrm{adapted}~\textrm{and}
    ~E\Big[\sup_{s\leq t\leq \tau}|X(t)|^2\Big]<+\infty\Big\},\\
&\mathscr{C}^2_{\bF}([s,\tau];H):=\Big\{X(\cdot):[s,\tau]\times\Omega\rrow
       H~|~\bF\textmd{-}\textrm{adapted with continuous path, and}~~~~~~~~~~\\
&~~~~~~~~~~~~~~~~~~~~~~~~~~~~~~~E\Big[\sup_{s\leq t\leq
\tau}|X(t)|^2\Big]<+\infty\Big\}.
\end{split}
\end{equation*}
For simplicity, define
$$\sH^2(s,\tau):=\sS^2_{\bF}([s,\tau];\bR^n)\times\sL^2_{\bF}(s,\tau;\bR^{n\times d})$$
equipped with norm
$$\|(\theta,\zeta)\|_{\sH^2(s,\tau)}:=\left\{E\Big[\sup_{s\leq t\leq
\tau}|\theta(t)|^2+\int_s^{\tau}|X(t)|^2\,dt\Big]\right\}^{\frac{1}{2}}.$$

Fix $\delta\geq0$ to be a constant. For all
$(y,z)\in\sH^2(0,T+\delta)$, let $(y(t),z(t))$ denote the value of
$(y,z)$ at time t, and $(y_t,z_t)$ denote the restriction of the
path of $(y,z)$ on $[t,t+\delta]$.

We give some conditions on the generator $(G,f)$ of NBSFDE
\eqref{equ}. Suppose that

$\bullet$ there exist two $\sB([0,T]\times
\bS^2([0,\delta];\bR^n)\times \bL^2 (0,\delta;\bR^{n\times
d}))\times \sF_T$-measurable functionals $J,F:
[0,T]\times\Omega\times \bS^2([0,\delta];\bR^n)\times
\bL^2(0,\delta;\bR^{n\times d})\rrow\bR^n$, such that
$J(\cdot,\vf,\psi)$ and $F(\cdot,\vf,\psi)$ are $\bF$-progressively
measurable, $\forall (\vf,\psi)\in \bS^2([0,\delta];\bR^n)\times
\bL^2(0,\delta;\bR^{n\times d})$.

$\bullet$ for all $(y,z)\in\sH^2(0,T+\delta)$,
$$G(t,y_t,z_t):=E\big[J(t,y_t,z_t)|\sF_t\big]=E_t\big[J(t,y_t,z_t)\big],$$
$$f(t,y_t,z_t):=E\big[F(t,y_t,z_t)|\sF_t\big]= E_t\big[F(t,y_t,z_t)\big].$$

Consider the following standing assumptions on $(G,f)$.

$(H1)$ There exist $\kappa\in(0,1)$ and a probability measure
$\lambda_1$ on $[0,\delta]$, such that for all
$(y,z),(\bar{y},\bar{z})\in\sH^2(0,T+\delta)$ and all
$(t,\omega)\in[0,T]\times\Omega$,
\begin{equation}\label{G}
\begin{split}
&~~|G(t,y_t,z_t)-G(t,\bar{y}_t,\bar{z}_t)|^2        \\
\leq\,&\kappa\,
E_t\Big[\int_0^{\delta}|y(t+u)-\bar{y}(t+u)|^2\,\lambda_1(du)
+\int_0^{\delta}|z(t+u)-\bar{z}(t+u)|^2\,du\Big].
\end{split}
\end{equation}

$(H2)$ There exist $L>0$ and two probability measures
$\lambda_2,~\lambda_3$ on $[0,\delta]$, such that for all
$(y,z),(\bar{y},\bar{z})\in\sH^2(0,T+\delta)$ and all
$(t,\omega)\in[0,T]\times\Omega$,
\begin{equation}\label{f}
\begin{split}
&~~~|f(t,y_t,z_t)-f(t,\bar{y}_t,\bar{z}_t)|^2   \\
\leq\,&
L\,E_t\left[\int_0^{\delta}|y(t+u)-\bar{y}(t+u)|^2\lambda_2(du)
      +\int_0^{\delta}|z(t+u)-\bar{z}(t+u)|^2\lambda_3(du)\right].
\end{split}
\end{equation}
$$(H3)~E\Bigl[\sup_{0\leq t\leq
T}|G(t,0,0)|^2\Bigr]<+\infty~~~~~\textrm{and}~~~~~~~E\Big[\int_0^T|f(t,0,0)|^2dt\Big]<+\infty,~~~~~~~~~~~~~~~~~~$$
where 0 is the path of $0\in\bS^2([0,\delta];\bR^n)$ or
$0\in\bL^2(0,\delta;\bR^n).$

\begin{rmk}
The dependence of G and f on the path of $z$ requirs in \eqref{f}
and \eqref{G} in different ways. In particular, such a form
$G(t,y_t,z(t))$ is not included in \eqref{G}.

Probability measures $\lambda_i~(i=1,2,3)$ in the generator allow us
to incorporate many interesting cases, such as
$$G(t,y_t,z_t)=\frac{1}{2\delta}E_t\Bigl[\int_0^{\delta}y(t+u)\,du\Bigr],$$
$$f(t,y_t,z_t)=E_t\Bigl[\int_0^{\delta}y(t+u)\,\lambda(du)+\int_0^{\delta}z(t+u)\,\lambda(du)\Bigr],$$
and
$$G(t,y_t,z_t)=E_t\Big[g(t,y(t),y(t+\theta_1),\int_0^{\delta}z(t+u)du)\Big],$$
$$f(t,y_t,z_t)=E_t\big[h(t,y(t),z(t),y(t+\theta_2),z(t+\theta_3))\big],$$
where $\theta_i>0$, $i=1,2,3$ are constants.
\end{rmk}

Let us introduce the definition of an adapted solution to NBSFDE
\eqref{equ}:
\begin{defn}
A pair of processes $(Y,Z)\in\sH^2(0,T+\delta)$ is called an adapted
$\bL^2$-solution to NBSFDE \eqref{equ} if they satisfy \eqref{equ}
in It\^{o}'s sense.
\end{defn}

The following theorem is devoted to the existence and uniqueness of
the adapted $\bL^2$-solution to NBSFDE \eqref{equ}.

\begin{thm}\label{ex_uni}
Let $(H1)$-$(H3)$ hold. Then for any pair
$(\xi,\zeta)\in\sH^2(T,T+\delta)$, NBSFDE \eqref{equ} admits a
unique adapted $\bL^2$-solution $(Y,Z)\in\sH^2(0,T+\delta)$.

Moreover, the following estimate holds:
\begin{equation*}
\begin{split}
&~~~E\left[\sup_{0\leq t\leq T}|Y(t)|^2+\int_0^T\!\!|Z(t)|^2\,dt\right]\\
\leq\,&C\,E\left[\sup_{T\leq t\leq T+\delta}|\xi(t)|^2
+\int_T^{T+\delta}\!\!\!|\zeta(t)|^2\,ds +\sup_{0\leq t\leq
T}|G(t,0,0)|^2+\int_0^T\!\!|f(t,0,0)|^2\,ds\right],
\end{split}
\end{equation*}
where $C$ only depends on $\,n,d,L$, and $\kappa$.
\end{thm}

\begin{proof}

\textbf{Step 1}. Define a subset of $\sH^2(0,T+\delta)$,
$$\sH^2(0,T;\xi,\zeta):=\left\{(y,z)\in\sH^2(0,T+\delta)
~\Bigm|~y(t)=\xi(t),~z(t)=\zeta(t), ~\forall
t\in[T,T+\delta]\right\}$$ equipped with norm
$$\|(y,z)\|^2=E\left[\sup_{0\leq t\leq T} e^{\beta t}|y(t)|^2+\int_0^T e^{\beta t}|z(t)|^2\,dt\right],$$
where $\beta>0$ is a constant waiting to be determined. It is
obvious that $\sH^2(0,T;\xi,\zeta)$ is a closed subset of
$\sH^2(0,T+\delta)$.

For $(y,z)\in\sH^2(0,T;\xi,\zeta)$, consider the equation
\begin{equation}\label{equ1}
\left\{
\begin{split}
&-d\big[Y(t)-G(t,y_t,z_t)\big]=f(t,y_t,z_t)\,dt-Z(t)\,dW_t,~~t\in[0,T];\\
&Y(t)=\xi(t),~~~~~~t\in[T,T+\delta];\\
&Z(t)=\zeta(t),~~~~~~t\in[T,T+\delta].
\end{split}
\right.
\end{equation}
Denote $\bar{Y}(t):=Y(t)-G(t,y_t,z_t)$, then
\begin{equation}\label{1bsde}
\bar{Y}(t)=\xi(T)-G(T,y_T,z_T)+\int_t^Tf(s,y_s,z_s)\,ds-\int_t^TZ(s)\,dW_s,~~t\in[0,T].
\end{equation}
From $(H1)$-$(H3)$, we have
\begin{equation*}
\begin{split}
&~~~E\left[|\xi(T)-G(T,y_T,z_T)|^2\right]\\
\leq\,&C\,E\left[|\xi(T)|^2+|G(T,0,0)|^2
         +\int_0^{\delta}|\xi(T+u)|^2\,\lambda_1(du)+\int_0^{\delta}|\zeta(T+u)|^2\,du\right]\\
\leq\,&C\,E\left[\sup_{T\leq t\leq
         T+\delta}|\xi(t)|^2+\int_T^{T+\delta}|\zeta(t)|^2\,dt\right]< +\infty
\end{split}
\end{equation*}
and
\begin{equation*}
\begin{split}
   &~~~E\left[\int_0^T|f(t,y_t,z_t)|^2\,dt\right]\\
\leq\,&C\,E\left[\int_0^T\left(|f(t,0,0)|^2+\int_0^{\delta}\!\!\!|y(t\!+\!u)|^2\lambda_2(du)
      +\int_0^{\delta}\!\!\!|z(t\!+\!u)|^2\lambda_3(du)\right)\,dt\right]\\
\leq\,&C\,E\left[\int_0^T|f(t,0,0)|^2\,dt+\sup_{0\leq t\leq
T+\delta}|y(t)|+\int_0^{T+\delta}|z(t)|^2\,dt\right]\,dt
         <+\infty.
\end{split}
\end{equation*}
In view of the theory of BSDEs, Eq.\eqref{1bsde} admits a unique
solution $(\bar{Y},Z)\in\sH^2(0,T)$. Then define
\begin{equation*}
Y(t):=\left\{
\begin{split}
&\bar{Y}(t)-G(t,y_t,z_t),~~~~~t\in[0,T];\\
&\xi(t),~~~~~~~~~~t\in[T,T+\delta],
\end{split}
\right.
\end{equation*}
and $Z(t):=\zeta(t),~t\in[T,T+\delta]$.

From $(H1)$ and $(H3)$ on G, we have
$$E\left[\sup_{0\leq t\leq T}|Y(t)|^2\right]\leq CE\left[\sup_{0\leq t\leq
T}|\bar{Y}(t)|^2+\sup_{0\leq t\leq
T+\delta}|y(t)|+\int_0^{T+\delta}|z(t)|^2\,dt\right]<\infty.$$ Then
$(Y,Z)\in\sH^2(0,T;\xi,\zeta)$ is the solution of Eq.\eqref{equ1}.

\textbf{Step 2}. Define a map $\Psi$ from $\sH^2(0,T;\xi,\zeta)$
onto itself. That is, $$\Psi:(y,z)\mapsto(Y,Z)$$ with $(Y,Z)$ being
the solution of \eqref{equ1} in Step 1. We prove $\Psi$ is a
contraction.

Take another pair of process $(\bar{y},\bar{z})\in
\sH^2(0,T;\xi,\zeta)$, and denote
$(\bar{Y},\bar{Z}):=\Psi(\bar{y},\bar{z})$. Let $\Delta Y(t):=
Y(t)-\bar{Y}(t)$, $\Delta Z(t):= Z(t)-\bar{Z}(t)$, $\Delta y(t):=
y(t)-\bar{y}(t),~\Delta z(t):= z(t)-\bar{z}(t)$ and $\Delta G(t):=
G(t,y_t,z_t)-G(t,\bar{y}_t,\bar{z}_t)$. Then
$$\Delta Y(t)-\Delta G(t)
=\int_t^T[f(s,y_s,z_s)-f(s,\bar{y}_s,\bar{z}_s)]\,ds+\int_t^T\Delta
Z(t)\,dW(t),~t\in[0,T].$$

Applying the It\^{o}'s formula for $e^{\beta t}|\Delta Y(t)-\Delta
G(t)|^2$, we have
\begin{equation}\label{Ito}
\begin{split}
 &e^{\beta t}|\Delta Y(t)-\Delta G(t)|^2
      +\int_t^T\beta e^{\beta s}|\Delta Y(s)-\Delta G(s)|^2\,ds+\int_t^Te^{\beta s}|\Delta Z(s)|^2\,ds\\
=\,&2\int_t^Te^{\beta s}\big\la\Delta Y(s)-\Delta
G(s),\,f(s,y_s,z_s)
      -f(s,\bar{y}_s,\bar{z}_s)\big\ra\,dt\\
&~~~-2\int_t^Te^{\beta s}
      \big\la\Delta Y(s)-\Delta G(s),\,\Delta Z(s)\big\ra\,dW(s).
\end{split}
\end{equation}
In view of $(H2)$ and the Schwartz inequality,
\begin{equation*}
\begin{split}
 &E\left[\int_0^T\beta e^{\beta s}
       |\Delta Y(s)-\Delta G(s)|^2\,ds+\int_0^Te^{\beta s}|\Delta Z(s)|^2\,ds\right]\\
=\,&2\,E\left[\int_0^Te^{\beta s}\big\la\Delta Y(s)-\Delta
   G(s),\,f(s,y_s,z_s)-f(s,\bar{y}_s,\bar{z}_s)\big\ra\,ds\right]\\
\leq& E\bigg[\int_0^TLCe^{\beta s}|\Delta Y(s)-\Delta G(s)|^2\,ds+
    \frac{1}{C}\int_0^Te^{\beta s}\Big(\int_0^{\delta}\!\!\!|\Delta
    y(s\!+\!u)|^2\lambda_2(du)\\
&~~~ +\int_0^{\delta}\!\!\!|\Delta
z(s\!+\!u)|^2\lambda_3(du)\Big)\,ds\bigg].
\end{split}
\end{equation*}
Consequently, choosing $\beta>LC$, we have
\begin{equation}\label{Ito1}
\begin{split}
      &~~~~E\left[\int_0^Te^{\beta s}|\Delta Z(s)|^2\,ds\right]\\
\leq\,&\frac{1}{C}\,E\left[\int_0^T e^{\beta
        s}\left(\int_0^{\delta}\!\!\!|\Delta y(s\!+\!u)|^2\lambda_2(du)
        +\int_0^{\delta}\!\!\!|\Delta
        z(s\!+\!u)|^2\lambda_3(du)\right)\,ds\right].
\end{split}
\end{equation}

From \eqref{Ito}, we have
\begin{equation*}
\begin{split}
      &E\left[\sup_{0\leq t\leq T}e^{\beta t}|\Delta Y(t)\!-\!\Delta G(t)|^2
         +\!\int_0^T\!\!\beta e^{\beta s}|\Delta Y(s)-\Delta G(s)|^2ds
         +\!\int_0^T\!\!e^{\beta s}|\Delta Z(s)|^2ds\right]\\
\leq\,&2\,E\left[\int_0^Te^{\beta s}|\Delta Y(s)-\Delta
          G(s)|\,|f(s,y_s,z_s)-f(s,\bar{y}_s,\bar{z}_s)|\,ds\right]\\
      &~~~~~~~~~~+2\,E\left[\sup_{0\leq t\leq T}|\int_t^Te^{\beta s}
            (\Delta Y(s)-\Delta G(s))\Delta Z(s)\,dW(s)|\right]\\
\leq\,&E\bigg[LC\,\int_0^Te^{\beta s}|\Delta Y(s)-\Delta G(s)|^2\,ds\\
      &~~~~~~~~~~+\frac{1}{C}\int_0^Te^{\beta s}\left(\int_0^{\delta}\!\!\!|\Delta
            y(s\!+\!u)|^2\lambda_2(du)
            +\int_0^{\delta}\!\!\!|\Delta z(s\!+\!u)|^2\lambda_3(du)\right)\,ds\bigg]\\
&~~~~~~~~~~+E\left[\frac{1}{K}\sup_{0\leq t\leq T}e^{\beta t}|\Delta
            Y(t)-\Delta G(t)|^2 +K\,\int_0^Te^{\beta s}|\Delta
            Z(s)|^2\,ds\right].
\end{split}
\end{equation*}
In view of \eqref{Ito1},
\begin{equation*}
\begin{split}
&(1-\frac{1}{K})E\left[\sup_{0\leq t\leq T}e^{\beta t}|\Delta
Y(t)-\Delta G(t)|^2\right]
+E\left[\int_0^Te^{\beta s}|\Delta Z(s)|^2ds\right]\\
\leq&\frac{K+1}{C}E\left[
    \int_0^T e^{\beta s}\left(\int_0^{\delta}\!\!\!|\Delta y(u\!+\!s)|^2\lambda_2(du)
+\int_0^{\delta}\!\!\!|\Delta
z(u\!+\!s)|^2\lambda_3(du)\right)\,ds\right].
\end{split}
\end{equation*}
Then we obtain
\begin{equation}\label{yz}
\begin{split}
&E\left[\sup_{0\leq t\leq T}e^{\beta t}|\Delta Y(t)-\Delta
G(t)|^2\right]
+E\left[\int_0^Te^{\beta s}|\Delta Z(s)|^2ds\right]\\
\leq&\theta E\left[
    \int_0^T e^{\beta s}\left(\int_0^{\delta}\!\!\!|\Delta y(s\!+\!u)|^2\lambda_2(du)
+\int_0^{\delta}\!\!\!|\Delta
z(s\!+\!u)|^2\lambda_3(du)\right)\,ds\right],
\end{split}
\end{equation}
where $\theta=\frac{K+1}{K-1}\frac{K}{C}$ can be any positive
constant by a proper choice of $K$ and $C$.

Since for all $a,~b\in\bR^n$ and $\forall \alpha\in(0,1)$,
$$(a-b)^2\geq (|a|-|b|)^2\geq
(1-\alpha)a^2-(\frac{1}{\alpha}-1)b^2,$$ then
\begin{equation*}
\begin{split}
      &\sup_{0\leq t\leq T}e^{\beta t}|\Delta Y(t)-\Delta G(t)|^2\\
\geq\,&\sup_{0\leq t\leq T}e^{\beta t}[(1-\alpha)|\Delta
     Y(t)|^2-(\frac{1}{\alpha}-1)|\Delta G(t)|^2]\\
\geq\,&(1-\alpha)\sup_{0\leq t\leq T}e^{\beta t}|\Delta
     Y(t)|^2-(\frac{1}{\alpha}-1)\sup_{0\leq t\leq T}e^{\beta t}|\Delta G(t)|^2].\\
\end{split}
\end{equation*}
\eqref{yz} becomes
\begin{equation*}
\begin{split}
&(1-\alpha)\,E\left[\sup_{0\leq t\leq T}e^{\beta
     t}|\Delta Y(t)|^2\right]+E\left[\int_0^Te^{\beta s}|\Delta Z(s)|^2\,ds\right]\\
\leq\,&(\frac{1}{\alpha}-1)\,E\left[\sup_{0\leq t\leq
           T}e^{\beta t}|G(t,y_t,z_t)-G(t,\bar{y}_t,\bar{z}_t)|^2\right]\\
&~~+\theta\, E\left[\int_0^T e^{\beta
          s}\left(\int_0^{\delta}\!\!|\Delta y(s\!+\!u)|^2\lambda_2(du)
          +\int_0^{\delta}\!\!\!|\Delta z(s\!+\!u)|^2\lambda_3(du)\right)\,ds\right]\\
\leq\,&(\frac{1}{\alpha}-1)\kappa\,E\left[\sup_{0\leq
      t\leq T}e^{\beta t}\left(\int_0^{\delta}\!\!|\Delta y(t\!+\!u)|^2\lambda_1(du)
      +\int_0^{\delta}\!\!|\Delta z(t\!+\!u)|^2du\right)\right]\\
&~~+\theta\,E\left[\int_0^T e^{\beta
     s}\left(\int_0^{\delta}\!\!\!|\Delta y(s\!+\!u)|^2\lambda_2(du)
         +\int_0^{\delta}\!\!\!|\Delta z(s\!+\!u)|^2\lambda_3(du)\right)\,ds\right]\\
\leq&[(\frac{1}{\alpha}-1)\kappa+\theta T]
        \,E\left[\sup_{0\leq t\leq T}e^{\beta t}|\Delta y(t)|^2\right]+[(\frac{1}{\alpha}-1)\kappa
      +\theta]\,E\left[\int_0^Te^{\beta t}|\Delta
      z(t)|^2\,dt\right].
\end{split}
\end{equation*}

To show that $\Psi$ is a contraction, it suffices to prove: for any
$\kappa\in(0,1)$, $\exists \alpha\in(0,1)$, $\theta>0$, such that
\begin{equation}\label{contact}
(\frac{1}{\alpha}-1)\kappa+\theta T<1-\alpha~~~~~~ \textrm{and}
~~~~~~(\frac{1}{\alpha}-1)\kappa+\theta<1. \end{equation} Indeed,
for any $\kappa\in(0,1)$, choose $\alpha\in(\kappa,1)$ and $\theta$
small sufficiently, the above two inequalities is easy to hold.

Therefore, $\Psi$ admits a unique fixed point. That is, \eqref{equ1}
admits a unique solution $(Y,Z)\in\sH^2(0,T;\xi,\zeta)$. In view of
the definition of $\sH^2(0,T;\xi,\zeta)$,
$(Y,Z)\in\sH^2(0,T+\delta)$ and it is the unique adapted
$\bL^2$-solution of NBSFDE \eqref{equ}.

\textbf{Step 3}. The estimation.

Let $(Y,Z)\in\sH^2(0,T+\delta)$ be the solution of Eq.\eqref{equ},
then
\begin{equation*}
\begin{split}
&e^{\beta t}|Y(t)-G(t,Y_t,Z_t)|^2+\int_t^T\!\!\beta e^{\beta
       s}|Y(s)-G(s,Y_t,Z_t)|^2\,ds+\int_t^T\!\!e^{\beta s}|Z(s)|^2\,ds\\
=\,&e^{\beta T}|\xi(T)-G(T,Y_T,Z_T)|^2+2\int_t^Te^{\beta
      s}\,\big\la Y(s)\!-\!G(s,Y_s,Z_s),\,f(s,Y_s,Z_s)\big\ra\,ds\\
&~~~~~~~-2\int_t^T\!\!e^{\beta
      s}\,\big\la Y(s)\!-\!G(s,Y_s,Z_s),\,Z(s)\big\ra\,dW(s)
\end{split}
\end{equation*}
Similar to the method in step 2, for all $\alpha\in(0,1)$ and $M>0$,
\begin{equation*}
\begin{split}
&(1-\alpha)E\Big[\sup_{0\leq t\leq T}e^{\beta
        t}|Y(t)|^2\Big]+E\Big[\int_0^T\!e^{\beta s}|Z(s)|^2\,ds\Big]\\
\leq&\,E\Big[e^{\beta T}|\xi(T)-G(T,Y_T,Z_T)|^2+\int_0^T\!e^{\beta
        t}|f(t,0,0)|^2\,dt\Big]+(\frac{1}{\alpha}\!-\!1)E\Big[\!\sup_{0\leq
        t\leq T}\!\!e^{\beta t}|G(t,Y_t,Z_t)|^2\Big]\\
&~~+\theta E\Big[\int_0^T\! e^{\beta
s}\Big(\int_0^{\delta}\!\!\!|Y(s\!+\!u)|^2\lambda_2(du)
        +\int_0^{\delta}\!\!\!|Z(s\!+\!u)|^2\lambda_3(du)\Big)\,ds\Big]\\
\leq&\,C\,E\Big[e^{\beta T}\Big(\!|G(T,0,0)|^2
        +\!\!\sup_{T\leq t\leq T+\delta}\!\!|\xi(t)|^2
       +\int_T^{T+\delta}\!\!|\zeta(t)|^2dt\!\Big)\\
&~~ +\int_0^T\!\!e^{\beta t}|f(t,0,0)|^2dt\Big]
        +(1\!+\!M)(\frac{1}{\alpha}\!-\!1)E\Big[\sup_{0\leq
        t\leq T}e^{\beta t}|G(t,0,0)|^2\Big]\\
&~~+(1\!+\!\frac{1}{M})(\frac{1}{\alpha}\!-\!1)\kappa
E\Big[\sup_{0\leq
        t\leq T}e^{\beta t}\Big(\int_0^{\delta}\!\!\!|Y(t\!+\!u)|^2\lambda_1(du)+\int_0^{\delta}\!\!\!|Z(t\!+\!u)|^2du\Big)\Big]   \\
&~~+\theta E\Big[\int_0^T \!\!e^{\beta
            s}\Big(\int_0^{\delta}\!\!\!|Y(u\!+\!s)|^2\lambda_2(du)
        +\int_0^{\delta}\!\!\!|Z(u\!+\!s)|^2\lambda_3(du)\Big)\,ds\Big]\\
\leq&\,C\,E\Big[e^{\beta T}\Big(\!\!\sup_{T\leq t\leq
T+\delta}\!\!|\xi(t)|^2
       +\!\!\int_T^{T+\delta}\!\!|\zeta(t)|^2dt\!\Big)+\!\!\int_0^T\!\!e^{\beta t}|f(t,0,0)|^2dt
       +\!\!\sup_{0\leq t\leq T}e^{\beta t}|G(t,0,0)|^2\Big]\\
&~+\bigl[(1\!+\!\frac{1}{M})(\frac{1}{\alpha}\!-\!1)\kappa\!+\!\theta
T\bigr]E\Big[\sup_{0\leq
        t\leq T+\delta}\!\!e^{\beta t}|Y(t)|^2\Big]
        +\bigl[(1\!+\!\frac{1}{M})(\frac{1}{\alpha}\!-\!1)\kappa\!+\!\theta\bigr]
        E\Big[\!\int_0^{T+\delta}\!\!e^{\beta t}|Z(t)|^2dt\Big].\\
\end{split}
\end{equation*}
Then,
\begin{equation*}
\begin{split}
&~~\bigl[(1\!-\!\alpha)\!-\!(1\!+\!\frac{1}{M})(\frac{1}{\alpha}\!-\!1)\kappa\!-\!\theta
T\bigr]E\Big[\sup_{0\leq t\leq T}e^{\beta t}|Y(t)|^2\Big]\\
&~~~~+\bigl[1\!-\!(1\!+\!\frac{1}{M})(\frac{1}{\alpha}\!-\!1)\kappa\!-\!\theta]
        E\Big[\int_0^Te^{\beta s}|Z(s)|^2ds\Big]\\
\leq&\,C\,E\Big[\sup_{T\leq T\leq
       T+\delta}\!\!|\xi(t)|^2+\int_T^{T+\delta}\!\!\!|\zeta(t)|^2dt+\int_0^T
       \!\!e^{\beta t}|f(t,0,0)|^2dt+\sup_{0\leq
        t\leq T}\!\!|G(t,0,0)|^2\Big].
\end{split}
\end{equation*}
Similar to the argument in step 2, for any $\kappa\in(0,1)$, there
exist $\alpha\in(0,1)$ and $M>0$, such that
$$(1-\alpha)-(1+\frac{1}{M})(\frac{1}{\alpha}-1)\kappa-\theta T>0;$$
$$1-(1+\frac{1}{M})(\frac{1}{\alpha}-1)\kappa-\theta>0.$$
Therefore,
\begin{equation*}
\begin{split}
      &E\Big[\sup_{0\leq t\leq T}|Y(t)|^2+\int_0^T|Z(s)|^2ds\Big]\\
\leq\,&C\,E\Big[\sup_{T\leq T\leq T+\delta}\!\!|\xi(t)|^2
       +\int_T^{T+\delta}\!\!\!|\zeta(t)|^2dt+\int_0^T\!\!|f(t,0,0)|^2dt
       +\sup_{0\leq t\leq T}\!\!|G(t,0,0)|^2\Big].
\end{split}
\end{equation*}
\end{proof}

\begin{rmk}
In \eqref{G}, $\kappa\in(0,1)$ is essential for the existence and
uniqueness. It is difficult to prove that $\Psi$ is a contraction
for $\kappa\geq1$, since \eqref{contact} does not hold for any
$\alpha\in(0,1)$ in this case,
\end{rmk}

\begin{rmk}\label{Z_terminal}
The value of $Z$ on $[0,T]$ is endogenous, just like BSDEs. The
terminal condition of $Z$ is needed for the well-posedness only on
$[T,T+\delta^\prime]$, where $\delta^\prime$ is the forward length
of the anticipation on Z of the generator $(G,f)$.

For example, if $G(t,\cdot,\cdot)$ and $f(t,\cdot,\cdot)$ only
depend on $Z(t)$, then the terminal condition of $Z$ is not
necessary. In this case, the solution $Z$ only makes sense on
$[0,T]$. But for the uniqueness of the solution in
$\sH^2(0,T+\delta)$, we define $Z=0$ on $[T,T+\delta]$, and omit it
in the equation for simplicity.
\end{rmk}

\begin{rmk}
(i) Since $Y(t)+G(t,Y_t,Z_t)$ is path-continuous, the continuity of
$Y(\cdot)$ depends on $G(t,\cdot,\cdot)$. If $G(t,\cdot,\cdot)$ is
continuous in t, for example,
$$G(t,Y_t,Z_t)=\frac{1}{2\delta}E_t\left[\int_t^{t+\delta}Y(s)ds\right],$$ then
$Y(\cdot)$ is also path-continuous.

(ii) $Y(t)+G(t,Y_t,Z_t)$ is a semi-martingale with diffusion
$\int_0^tZ(s)\,dW(s)$. However, we do not know whether $Y(\cdot)$ is
semi-martingale or not.
\end{rmk}

Similar to the proof of the former theorem, we have the following
corollary:
\begin{cor}\label{depend_gener}
Suppose that $(G_i,f_i)$ satisfies $(H1)$-$(H3)$, and
$(\xi^i,\zeta^i)\in\sH(T,T+\delta),~i=1,2$. Consider the following
NBSFDE:
\begin{equation*}
\left\{
\begin{split}
&-d\,[Y(t)-G_i(t,Y_t,Z_t)]=f_i(t,Y_t,Z_t)\,dt-Z(t)\,dW_t,~~~~~t\in[0,T];\\
&Y(t)=\xi^i(t),~~~~~~t\in[T,T+\delta];\\
&Z(t)=\zeta^i(t),~~~~~~t\in[T,T+\delta].
\end{split}
\right.
\end{equation*}
Let $(Y^i,Z^i)\in\sH(0,T+\delta)$ be the solution for $i=1,2$,
respectively. Then the following estimate holds:
\begin{equation}\label{estim}
\begin{split}
&~~~E\left[\sup_{0\leq t\leq T}|Y^1(t)-Y^2(t)|^2+\int_0^T |Z^1(t)-Z^2(t)|^2dt\right]\\
\leq\,&C\,E\bigg[\sup_{T\leq t\leq T+\delta}|\xi^1(t)-\xi^2(t)|^2
       +\int_T^{T+\delta} |\zeta^1(t)-\zeta^2(t)|^2\,dt\\
&+\!\sup_{0\leq t\leq
T}|G_1(t,Y^1_t,Z^1_t)\!-\!G_2(t,Y^1_t,Z^1_t)|^2
        +\int_0^T\!\!|f_1(t,Y^1_t,Z^1_t)\!-\!f_2(t,Y^1_t,Z^1_t)|^2dt\bigg],
\end{split}
\end{equation}
where $C$ only depends on $\,n,d,L$ and $\kappa$.
\end{cor}

This corollary states the dependence of solutions on the generator
$(G,f)$ and the terminal conditions. In the following, we will
discuss how the solution depends on $\delta$ in a simple case.

Let $\delta_1$ and $\delta_2$ be two nonnegative constant, and
$\delta_1>\delta_2$. Consider the following two equations:
\begin{equation}\label{depend_delta}
\left\{\!\!
\begin{split}
&\!-d\Big[Y(t)-E_t\big[g(t,Y(t+\delta_i))\big]\!\Big]\!=\!E_t\big[h(t,Y(t),Z(t),Y(t+\delta_i))\big]dt-Z(t)\,dW(t);\\
&Y(t)=\xi(t),~~~~t\in[T,T+\delta_1].
\end{split}
\right.
\end{equation}
where $g:[0,T]\times\Omega\times\bR^n\rrow\bR^n$, and
$h:[0,T]\times\Omega\times\bR^n\times\bR^{n\times
d}\times\bR^n\rrow\bR^n$ are both adapted processes.

Since both g and h are independent of the anticipation of
$Z(\cdot)$, by Remark \ref{Z_terminal}, it is sufficient for the
well-posedness to give the terminal value of Y on $[T,T+\delta]$.

Similar to assumptions $(H1)$ and $(H2)$, consider the assumption

$(H4)$ There are $\kappa\in(0,1)$ and $L>0$, such that for all
$y,\bar{y},v,\bar{v}\in\bR^n$ and $z,\bar{z}\in\bR^{n\times d},$
$$|g(t,y)-g(t,\bar{y})|
\leq \kappa |y-\bar{y}|;$$
$$|h(t,y,z,v)-h(t,\bar{y},\bar{z},\bar{v})|\leq
L\left(|y-\bar{y}|+|z-\bar{z}|+|v -\bar{v}|\right).$$


Then we give the following proposition about the dependence of
solution on $\delta$:
\begin{prop}
Suppose that $(H4)$ holds, and assume there are $C>0$ and
$\alpha>0$, such that
$$|g(t,y)-g(t',y)|\leq C|t-t'|^{\alpha},~~ \forall y\in\bR^n,$$
and for any $\vf,\bar{\vf}\in\bL^2(\sF_T;\bR^n),$
$$|g(t,\vf)-g(t,\bar{\vf})|\leq \kappa |E_t[\vf-\bar{\vf}]|;$$
$$|h(t,y,z,\vf)-h(t,y,z,\bar{\vf})|\leq L |E_t[\vf-\bar{\vf}]|.$$

Then for any $\xi\in\sS^2_{\bF}([T,T+\delta_1];\bR^n)$, the
following estimate holds
\begin{equation*}
\begin{split}
&~~~~E\left[\sup_{0\leq t\leq T}|\Delta Y(t)|^2+\int_0^T|\Delta
          Z(t)|^2dt\right]\\
\leq&\, C\Big(|\delta_1-\delta_2|^{2\alpha}+|\delta_1-\delta_2|
+E\Big[\sup_{t,\bar{t}\in[T,T+\delta]}|E_T[\xi(t)-\xi(\bar{t})]|^2\Big]\Big).
\end{split}
\end{equation*}
Particularly, if $\xi$ is a martingale, then
$$E\left[\sup_{0\leq t\leq T}|\Delta Y(t)|^2+\int_0^T|\Delta Z(t)|^2dt\right]
\leq
C\left[|\delta_1-\delta_2|^{\alpha}+|\delta_1-\delta_2|\right].$$
\end{prop}

\begin{proof}
It is obvious that \eqref{depend_delta} admits a unique pair of
solution $(Y^i,Z^i)$ for $i=1,2$ respectively. 
Let $\Delta Y:= Y_1-Y_2$ and $\Delta Z:= Z_1-Z_2$. In view of
\eqref{estim}, we have the following estimate:
\begin{equation}\label{deltaesti}
\begin{split}
      &E\Big[\sup_{0\leq t\leq T}|\Delta Y(t)|^2+\int_0^T|\Delta Z(t)|^2dt\Big]\\
\leq\,&C\,E\Big[\sup_{0\leq t\leq T}|g(t,Y_1(t+\delta_1))-g(t,Y_1(t+\delta_2))|^2\\
       &~~~+\int_0^T|h(t,Y_1(t),Z_1(t),Y_1(t+\delta_1))-h(t,Y_1(t),Z_1(t),Y_1(t+\delta_2)|^2dt\Big]\\
\leq\,&C\,E\Big[\sup_{0\leq t\leq
       T}|E_t[Y_1(t+\delta_1)-Y_1(t+\delta_2)]|^2\Big].
\end{split}
\end{equation}

If $t+\delta_1\geq T,~t+\delta_2\geq T$, then
\begin{equation*}
\begin{split}
&\left|E_t\big[Y_1(t+\delta_1)-Y_1(t+\delta_2)\big]\right|^2
=\,\left|E_t\big[\xi(t+\delta_1)-\xi(t+\delta_2)\big]\right|^2 \\
\leq\,&
E_t\left|E_T\big[\xi(t+\delta_1)-\xi(t+\delta_2)\big]\right|^2\\
\leq\,&
E_t\Big[\sup_{s,\bar{s}\in[T,T+\delta]}|E_T[\xi(s)-\xi(\bar{s})]|^2\Big].
\end{split}
\end{equation*}

If  $t+\delta_1\geq T,~t+\delta_2<T$, then
$T-(t+\delta_2)\leq(t+\delta_1)-(t+\delta_2)=\delta_1-\delta_2,$
\begin{equation*}
\begin{split}
&E_t\big[Y_1(t+\delta_1)-Y_1(t+\delta_2)\big]=E_t\big[\xi(t+\delta_1)-Y_1(t+\delta_2)\big]\\
=\,&E_t\Big[\xi(t+\delta_1)-\xi(T)
+[g(T,\xi(T+\delta_1))-g(t+\delta_2,\xi(t+\delta_2+\delta_1))]\\
&~~~
-\int_{t+\delta_2}^Th(s,Y_1(s),Z_1(s),Y_1(s+\delta_1))\,ds\Big].
\end{split}
\end{equation*}
Then
\begin{equation*}
\begin{split}
&|E_t[Y_1(t+\delta_1)-Y_1(t+\delta_2)]|^2\\
\leq\,&C\,\Big(|E_t[\xi(t+\delta_1)-\xi(T)]|^2+|E_t[g(T,\xi(T+\delta_1))
           -g(t+\delta_2,\xi(t+\delta_2+\delta_1))]|^2\\
&~~~~+(T-t-\delta_2)E_t\Big[\int_{t+\delta_2}^T|h(s,Y_1(s),Z_1(s),Y_1(s+\delta_1))|^2ds\Big]\Big)\\
\leq\,&C\,\Big(E_t\Big[|E_T[\xi(t+\delta_1)-\xi(T)]|^2+|E_T[\xi(T+\delta_1)-\xi(t+\delta_1+\delta_2)]|^2\Big]\\
      &~~~~~~~~+|\delta_1-\delta_2|^{2\alpha}+|\delta_1-\delta_2|\Big)\\
\leq\,&C\,\Big(E_t\Big[\sup_{s,\bar{s}\in[T,T+\delta]}|E_T[\xi(s)-\xi(\bar{s})]|^2\Big]
         +|\delta_1-\delta_2|^{2\alpha}+|\delta_1-\delta_2|\Big).
\end{split}
\end{equation*}

If $t+\delta_1\leq T$, then
\begin{equation*}
\begin{split}
    &E_t[Y_1(t+\delta_1)-Y_1(t+\delta_2)]\\
=&\,E_t\left[g(t+\delta_1,Y_1(t+2\delta_1))-g(t+\delta_2,Y_1(t+\delta_2+\delta_1))\right]\\
&~~~~+E_t\Big[\int_{t+\delta_2}^{t+\delta_1}
      h(s,Y_1(s),Z_1(s),Y_1(s+\delta_1))ds\Big],
\end{split}
\end{equation*}
and
\begin{equation*}
\begin{split}
   &|E_t[Y_1(t+\delta_1)-Y_1(t+\delta_2)]|^2\\
\leq\,& \left|E_t\Big[g(t+\delta_1,Y_1(t+2\delta_1))-g(t+\delta_2,Y_1(t+\delta_2+\delta_1))\Big]\right|^2\\
&+|\delta_1-\delta_2|E_t\Big[\int_{t+\delta_2}^{t+\delta_1}
     |h(s,Y_1(s),Z_1(s),Y_1(s+\delta_1))|^2ds\Big]\\
\leq &\,C\Big(|\delta_1-\delta_2|^{2\alpha}
     +|E_t[Y_1(t+2\delta_1)-Y_1(t+\delta_1+\delta_2)]|^2+|\delta_1-\delta_2|\Big)\\
\leq\,&\cdots\,\leq C\Big(|\delta_1-\delta_2|^{2\alpha}
       +|\delta_1-\delta_2|+E_t\Big[\sup_{s,\bar{s}\in[T,T+\delta]}|E_T[\xi(s)-\xi(\bar{s})]|^2\Big]\Big).
\end{split}
\end{equation*}

So we obtain
\begin{equation*}
\begin{split}
&~~~E\left[\sup_{0\leq t\leq T}|E_t[Y_1(t+\delta_1)-Y_1(t+\delta_2)]|^2\right] \\
\leq&\, C\Big(|\delta_1-\delta_2|^{2\alpha}+|\delta_1-\delta_2|
        +E_t\Big[\sup_{s,\bar{s}\in[T,T+\delta]}|E_T[\xi(s)-\xi(\bar{s})]|^2\Big]\Big),~~~\forall
        t\in[0,T].
\end{split}
\end{equation*}

Applying \eqref{deltaesti}, then
\begin{equation*}
\begin{split}
    &~~~E\Big[\sup_{0\leq t\leq T}|\Delta Y(t)|^2+\int_0^T|\Delta
      Z(t)|^2dt\Big]\\
\leq&\, C\Big(|\delta_1-\delta_2|^{2\alpha}+|\delta_1-\delta_2|
+E_t\Big[\sup_{s,\bar{s}\in[T,T+\delta]}|E_T[\xi(s)-\xi(\bar{s})]|^2\Big]\Big).
\end{split}
\end{equation*}

If $\xi(\cdot)$ is a martingale, then $E_T[\xi(t)]=\xi(T)$,
$$E\Big[\sup_{0\leq t\leq
T}|\Delta Y(t)|^2+\int_0^T|\Delta Z(t)|^2dt\Big] \leq
C\big(|\delta_1-\delta_2|^{\alpha}+|\delta_1-\delta_2|\big). $$
\end{proof}

\section{Comparison Theorem}
In this section, we are going to establish a comparison theorem for
the adapted solution of NBSFDEs when $G$ is independent of Z and $f$
depends on $Z$ without anticipation, i.e.
\begin{equation}\label{equ2}
\left\{
\begin{split}
&-d\big[Y(t)-G(t,Y_t)\big]=f(t,Y_t,Z(t))\,dt-Z(t)\,dW(t),~~t\in[0,T];\\
&Y(t)=\xi(t),~~~~~~t\in[T,T+\delta].
\end{split}
\right.
\end{equation}
In view of Remark \ref{Z_terminal}, we omit the terminal condition
on $Z$, and treat $Z=0$ on $[T,T+\delta]$. Here we suppose that
$n=1,~d\geq1$.

Both assumptions $(H1)$ and $(H2)$ on $(G,f)$ are reduced to the
following:

$(H5)$ For any $(y,z),(\bar{y},\bar{z})\in\sH^2(0,T+\delta)$,
$$|G(t,y_t)-G(t,\bar{y}_t)|^2\leq
\kappa \,E_t\int_0^{\delta}|y(t+u)-\bar{y}(t+u)|^2\lambda_1(du),$$
$$|f(t,y_t,z(t))\!-\!f(t,\bar{y}_t,\bar{z}(t))|^2
\leq
L\,E_t\Bigl[\int_0^{\delta}\!\!|y(t+u)-\bar{y}(t+u)|^2\lambda_2(du)+|z(t)-\bar{z}(t)|^2\Big],$$
where $\kappa$, $L$, $\lambda_1$ and $\lambda_2$ are the same as in
$(H1)$ and $(H2)$.

Let us first consider the simple case:
\begin{equation}
\left\{
\begin{split}
&-d\,\big[y(t)-g(t)\big]=\big[f(t)+\beta(t)z(t)\big]\,dt-z(t)\,dW(t),~~t\in[0,T],\\
&Y(T)=Q\,\in\bL^2(\sF_T;\bR).
\end{split}
\right.
\end{equation}

We have the following lemma:
\begin{lem}\label{compar1}
Suppose that $Q\geq0$, $g(T)=0$. Then for any
$~\beta\in\sL^{\infty}_{\bF}(0,T;\bR)$, and nonnegative
$g\in\sS^2_{\bF}(0,T;\bR)$ and $f\in\sL^2_{\bF}(0,T;\bR)$, we have
almost surely $y(t)\geq0,~\forall t\in[0,T]$.
\end{lem}

\begin{proof}
Define $\overline{y}(t)=y(t)-g(t)$. Then
$$\overline{y}(t)=Q+\int_t^T\big[f(s)
+\beta(s)z(s)\big]\,ds-\int_t^Tz(s)\,dW(s).$$

Since $Q\geq0$ and $f(s)\geq0$, In view of the comparison theorem of
BSDEs, we have $\overline{y}(t)=y(t)-g(t)\geq0$, i.e. $y(t)\geq
g(t)\geq0$.
\end{proof}

Now, we have the following comparison theorem.
\begin{thm}\label{compar}
Suppose that $(H3)$ and $(H5)$ hold for $(f_i,G_i)$, and
$\xi^i\in\sS^2_{\bF}([T,T+\delta];\bR^n),~i=1,2$. Let $(Y^i,Z^i)$ be
the solution of the following NBSFDE, $i=1,2$ respectively:
\begin{equation}\label{com_equ}
\left\{
\begin{split}
&-d\,\big[Y^i(t)-G_i(t,Y^i_t)\big]=f_i(t,Y^i_t,Z^i(t))\,dt-Z^i(t)\,dW(t),~~t\leq T;\\
&Y^i(t)=\xi^i(t),~~~~~~t\in[T,T+\delta].
\end{split}
\right.
\end{equation}
Suppose that $G_i(T,\xi^j_T)=0$ $(i,j=1,2)$, and both $f_2$ and
$G_2$ are nondecreasing in the path of $Y$. If
$\xi^1(T)\geq\xi^2(T)$,
$$G_1(t,Y_t^1)\geq G_2(t,Y_t^1)~~\textrm{and}~~f_1(t,Y_t^1,Z^1(t))\geq f_2(t,Y_t^1,Z^1(t)),$$
then we have almost surely $Y^1(t)\geq Y^2(t)$ for all $t\in[0,T]$.
\end{thm}

\begin{proof}
In view of Step 1 in Theorem \ref{ex_uni}, the following equation
\begin{equation*}
\left\{
\begin{split}
&-d\,\big[Y^3(t)-G_2(t,Y_t^1)\big]=f_2(t,Y_t^1,Z^3(t))\,dt-Z^3(t)\,dW(t),~~~t\leq T;\\
&Y^3(t)=\xi^2(t),~~~~~~t\in[T,T+\delta]
\end{split}
\right.
\end{equation*}
admits a unique pair of solution $(Y^3,Z^3)$. Define $\Delta_3 Y:=
Y^1-Y^3,~\Delta_3 Z:= Z^1-Z^3,~\Delta G(t):=
G_1(t,Y_t^1)-G_2(t,Y_t^1)$,
$$\Delta f(t):=
f_1(t,Y_t^1,Z^1(t))-f_2(t,Y_t^1,Z^1(t))\geq0,$$ and
$$
\beta(t):=\frac{[f_2(t,Y_t^1,Z^1(t))
        -f_2(t,Y_t^1,Z^3(t))]\Delta_3Z(t)}{|\Delta_3Z(t)|^2}\bI_{\{|\Delta_3 Z(t)|\neq0\}}.$$
We have $\Delta G(t)\geq0,~\Delta f(t)\geq0$,
$\beta\in\sL^{\infty}_{\bF}(0,T;\bR^d)$, and
\begin{equation*}
\begin{split}
    &\Delta_3 Y(t)-\Delta G(t)\\
=\,&\Delta\xi(T)+\int_t^T[f_1(s,Y_s^1,Z^1(s))
        -f_2(s,Y_s^1,Z^3(s))]\,ds-\int_t^T\Delta_3 Z(s)dW(s)\\
=&\Delta \xi(T)+\int_t^T[\Delta f(s)+\beta^T(s)\Delta_3 Z(s)]\,ds
        -\int_t^T\Delta_3 Z(s)\,dW(s).
\end{split}
\end{equation*}
In view of Lemma \ref{compar1}, we have almost surely $\Delta_3
Y(t)=Y^1(t)-Y^3(t)\geq0, ~\forall t\in[0,T]$. Since both $G_2$ and
$f_2$ are nondecreasing in the path of $Y$, then
$$G_2(t,Y_t^1)\geq G_2(t,Y_t^3)~~\textrm{and}~~f_2(t,Y_t^1,\cdot)\geq
f_2(t,Y_t^3,\cdot).$$

The following equation
\begin{equation*}
\left\{
\begin{split}
&-d\,\big[Y^4(t)-G_2(t,Y_t^3)\big]=f_2(t,Y_t^3,Z^4(t))\,dt-Z^4(t)\,dW(t),~~t\leq T;\\
&Y^4(t)=\xi^2(t),~~~~~~t\in[T,T+\delta]
\end{split}
\right.
\end{equation*}
has a unique pair of solution $(Y^4,Z^4)$. Define $\Delta_4 Y:=
Y^3-Y^4,~\Delta_4 Z:= Z^3-Z^4$, then
\begin{equation*}
\begin{split}
&\Delta_4 Y(t)-[G_2(t,Y_t^1)-G_2(t,Y_t^3)]\\
=&\int_t^T\Big(f_2(s,Y_s^1,Z_3(s))-f_2(s,Y_s^3,Z_4(s))\Big)\,ds-\int_t^T\Delta_4Z(s)\,dW(s)\\
=&\!\!\int_t^T\!\Big(\!f_2(s,Y_s^1,Z_3(s))-f_2(s,Y_s^3,Z_3(s))+f_2(s,Y_s^3,Z_3(s))-f_2(s,Y_s^3,Z_4(s))\!\Big)\,ds\\
&~~~~~~~-\int_t^T\Delta_4Z(s)\,dW(s).
\end{split}
\end{equation*}
Similarly as in the preceding paragraph, we have almost surely
$\Delta_4Y(t)=Y_3(t)-Y_4(t)\geq0,~\forall t\in[0,T]$.

For any integer n, the following equation
\begin{equation}\label{sequ}
\left\{
\begin{split}
&-d\,\big[Y^n(t)-G_2(t,Y_t^{n-1})\big]=f_2(t,Y_t^{n-1},Z^n(t))\,dt-Z^n(t)\,dW(t),~~t\leq T;\\
&Y^n(t)=\xi^2(t),~~~~~~t\in[T,T+\delta]
\end{split}
\right.
\end{equation}
admits a unique pair of solution $(Y^n,Z^n)$. Similarly, we deduce
almost surely
$$Y^1(t)\geq Y^3(t)\geq Y^4(t)\geq\cdots\geq
Y^n(t)\geq\cdots,~\forall t\in[0,T].$$

In view of Steps 2 and 3 of Theorem \ref{ex_uni}, there exists
$\mu\in(0,1)$ such that
\begin{equation*}
\begin{split}
     &E\left[\sup_{0\leq t\leq T}e^{\beta t}|\Delta_nY(t)|^2
           +\int_0^Te^{\beta t}|\Delta_nZ(t)|^2\,dt\right]\\
\leq\ &\mu E\left[\sup_{0\leq t\leq T}e^{\beta
t}|\Delta_{n-1}Y(t)|^2
           +\int_0^Te^{\beta t}|\Delta_{n-1}Z(t)|^2\,dt\right]\\
\leq\ &\cdots\leq \mu^{n-2}E\left[\sup_{0\leq t\leq T}e^{\beta
            t}|\Delta_3Y(t)|^2+\int_0^T e^{\beta t}|\Delta_3Z(t)|^2\,dt\right]\\
\leq\ &\mu^{n-2}C\,E\left[|\xi^1-\xi^2|^2
       +\sup_{0\leq t\leq T}e^{\beta T}|\Delta G(t)|^2
       +\int_0^T e^{\beta t}|\Delta
       f(s)|^2dt\right]\longrightarrow\,0.
\end{split}
\end{equation*}
That is, $(Y^n,Z^n)_{n\geq1}$ is a Cauchy sequence in
$\sH^2(0,T+\delta)$.

Let $n$ go to infinity in Eq.\eqref{sequ} and compare it to Eq.
\eqref{com_equ}, we have $$(Y^n,Z^n)\rrow
(Y^2,Z^2),~\textrm{in}~\sH^2(0,T+\delta).$$ Then we have almost
surely $Y^1(t)\geq Y^2(t),~\forall t\in[0,T].$

\end{proof}

\section{Application in Optimal Stochastic Control}
As an application of the previous result, we discuss optimal control
of simple neutral stochastic functional differential equations
(NSFDEs) in this section.

Here we only discuss the optimal control of NSFDEs with a special
form as follow. The more general one will be discussed in the
future. Let $n=1$ and $d\geq1$ for simplicity,
\begin{equation}\label{NDSde}
\left\{\!\!
\begin{split}
&d\,\Big[X(t)-\kappa\int_0^{\delta}\!\!\!\!X(t\!-\!r)\mu_1(dr)\Big]=b(t,X(t),\int_0^{\delta}\!\!\!\!X(t\!-\!r)\mu_2(dr),u(t))\,dt\\
&~~~~~~~~~~+\Big[c(t,u(t))X(t)-\kappa\int_0^{\delta}\!\!\!c(t\!-\!r,u(t\!-\!r))X(t\!-\!r)\mu_1(dr)\Big]\,dW(t),~t\in[0,T];\\
&X(t)=\vf(t),~~~~~~t\in[-\delta,0],
\end{split}
\right.
\end{equation} where $\mu_1$ and $\mu_2$ are two probability measures in
$(0,\delta]$, $\kappa\in(-1,1)$ is a constant, and
$$b:[0,T]\times\Omega\times\bR\times\bR\times
\bR^m\rrow\bR,~\textrm{and }
c:[0,T]\times\Omega\times\bR^m\rrow\bR^d
$$ are jointly measurable, and
$b(\cdot,x,y,u)$ and $c(\cdot,u)$ are $\bF\textmd{-}$progressively
measurable for all $(x,y,u)\in\bR^n\times\bR^n\times\bR^m$.

The cost functional is defined as follow:
$$J(u(\cdot)):= E\left[\int_0^T l(t,X(t),u(t))\,dt+M(X(T))\right],$$
where $l:[0,T]\times\Omega\times\bR\times\bR^m\rrow\bR$ and
$M:\Omega\times\bR\rrow\bR$ are jointly measurable, and
$l(\cdot,x,u)$ is $\bF\textmd{-}$progressively measurable for all
$(x,u)\in\bR^n\times\bR^m$.

We introduce the following assumptions:

(H6) $b(t,\cdot,\cdot,\cdot),~c(t,\cdot),~l(t,\cdot,\cdot)$ and
$M(\cdot)$ are continuously differentiable in $(x,y,u)$ with bounded
derivatives.

Let $U\subseteq\bR^m$ be a nonempty bounded subset. Define the
admissible control set
$$\sU:=\Bigl\{u:[-\delta,T]\times\Omega\rrow U,~\bF\textmd{-}\textrm{progressively measurable}~|~u(t)=0,~\forall t\in[-\delta,0)\Bigr\}.$$
In view of the conclusions in \cite{Kolmanovskii1986},
Eq.\eqref{NDSde} admits an unique solution
$X(\cdot)\in\sS^2_{\bF}([-\delta,T];\bR)$ for all $u(\cdot)\in\sU$
and $\vf(\cdot)\in\bS([-\delta,0];\bR)$. Thus the cost functional is
well-defined.

Our optimal control problem can be stated as follow:

\textbf{Problem(C)}: Find a control process $\bar{u}(\cdot)\in\sU$
such that
$$J(\bar{u}(\cdot))=\sup_{u(\cdot)\in\sU}J(u(\cdot)).$$
Then $\bar{u}(\cdot)$ is called the optimal control, and
$(\bar{X}(\cdot),\bar{u}(\cdot))$ is called the optimal pair with
$\bar{X}(\cdot)$ being the solution of Eq.\eqref{NDSde}
corresponding to $\bar{u}(\cdot)$.

\subsection{The Pontryagin maximum principle}
In this subsection, we construct the Pontryagin maximum principle
for the former optimal control problem (C).

At first, we establish the duality between linear NSFDEs and linear
NBSFDEs, which is crucial in the solution of optimal control
problem.

For all $t\in[0,T]$, consider the following linear NSFDE:
\begin{equation}\label{ForwardDE}
\left\{\!
\begin{split}
&d\,\Big[X(s)-\kappa\int_0^{\delta}\!\!\!X(s\!-\!r)\mu_1(dr)\Big]=\Big[p(s)X(s)
      +q(s)\int_0^{\delta}\!\!\!X(s\!-\!r)\mu_2(dr)+f(s)\Big]\,ds\\
&~~~~~~~~+\Big[v(s)X(s)+\rho(s)-\kappa\int_0^{\delta}\!\!\!\big(v(s\!-\!r)X(s\!-\!r)\!
      +\!\rho(s\!-\!r)\big)\mu_1(dr)\Big]\,dW(s),~s\in[t,T];\\
&X(t)=x;\\
&X(s)=0,~~~~~s\in[t-\delta,t),
\end{split}
\right.
\end{equation}
and the linear NBSFDE
\begin{equation}\label{BFE}
\left\{
\begin{split}
&-d\Big[Y(s)-\kappa E_s\int_0^{\delta}\!\!\!Y(s\!+\!r)\mu_1(dr)\Big]
   =\Big[p(s)Y(s)+E_s\!\int_0^{\delta}\!\!\!q(s\!+\!r)Y(s\!+\!r)\mu_2(dr)\\
&~~~~~~~~~~~~~~~~~~~~~~~~~~~~~~~~~~~~~~~~~~~~~+v(s)Z(s)+h(s)\Big]\,ds-Z(s)\,dW(s),~~s\in[t,T];\\
&Y(T)=\xi;\\
&Y(s)=0,~~~~s\in(T,T+\delta],
\end{split}
\right.
\end{equation}
where $x\in\bR$, $\xi\in\bL^2(\sF_T;\bR)$,
$p,q\in\sL^{\infty}_{\bF}(0,T;\bR)$,
$v\in\sL^{\infty}_{\bF}(-\delta,T;\bR^d)$,
$f,h\in\sL^2_{\bF}(0,T;\bR)$, and
$\rho\in\sL^2_{\bF}(-\delta,T;\bR)$.

Then we have the following lemma.
\begin{lem}\label{duality}
Suppose that $\rho(s)=0$ for all $s<t$. Let $X$ and $Y$ be solution
of Eq.\eqref{ForwardDE} and Eq.\eqref{BFE}, respectively. Then
$$E_t\left[\xi X(T)+\int_t^T
h(s)X(s)\,ds\right]=xY(t)+E_t\left[\int_t^T\!\!
\Big(f(s)Y(s)+\rho(s)Z(s)\Big)\,ds\right].$$
\end{lem}

\begin{proof}
\begin{equation}\label{1}
\begin{split}
   &E_t\int_t^Td\, \Big\{\Big[X(s)-\kappa \int_0^{\delta}\!\!\!X(s\!-\!r)\mu_1(dr)\Big]Y(s)\Big\}\\
=\,&E_t\bigg\{\int_t^T\!Y(s)\,d\Big[X(s)\!-\!\kappa
\int_0^{\delta}\!\!\!X(s\!-\!r)\mu_1(dr)\Big]\!+\!\int_t^T\!\Big[X(s)\!-\!\kappa\!
\int_0^{\delta}\!\!\!X(s\!-\!r)\mu_1(dr)\Big]dY(s)\\
&~~~~ +\int_t^Td\Big[X(s)-\kappa \int_0^{\delta}\!\!\!X(s\!-\!r)\mu_1(dr)\Big]dY(s)\bigg\}\\
=\,&E_t\left[\int_t^T\Big(\big(
       p(s)X(s)+q(s)\int_0^{\delta}\!\!\!X(s\!-\!r)\mu_2(dr)\big)Y(s)+f(s)Y(s)\Big)\,ds\right]\\
   &+E_t\left[\int_t^TX(s)\,dY(s)
   -\int_t^T\kappa\int_0^{\delta}\!\!\!X(s\!-\!r)\mu_1(dr)\,dY(s)\right]\\
   &+E_t\left[\int_t^T\!\!\Big(v(s)X(s)\!+\!\rho(s)\!
      -\!\kappa\!\!\int_0^{\delta}\!\!\!\big(v(s\!-\!r)X(s\!-\!r)\!+\!\rho(s\!-\!r)\big)\mu_1(dr)\Big)dW(s)\,dY(s)\right]\\
=\,&E_t\bigg\{\int_t^TX(t)\Bigl[q(s)Y(s)+E_s\int_0^{\delta}\!\!\!q(s\!+\!r)Y(s\!+\!r)\mu_2(dr)\Bigr]\,ds\\
       &+\int_t^TX(s)\,d\,\Big[Y(s)-\kappa\int_0^{\delta}\!\!\!Y(s\!+\!r)\mu_1(dr)\Big]
       +\int_{T-\delta}^TX(s)\,d\,\Bigl[\kappa \int_0^{\delta}\!\!\!Y(s\!+\!r)\mu_1(du)\Bigr]\\
   &+\int_t^T\Big(v(s)X(s)\!+\!\rho(s)\Big)dW(s)\,d\Big[Y(s)-\kappa\int_0^{\delta}\!\!\!Y(s\!+\!r)\mu_1(dr)\Big]
   \!+\!\int_t^T\!\!f(s)Y(s)\,ds\bigg\}\\
=\,&E_t\Big[-\int_t^TX(s)h(s)\,ds-\kappa\int_0^{\delta}\!\!\!X(T\!-\!r)\mu_1(dr)\,\xi
       +\int_t^T\Big(f(s)Y(s)+\rho(s)Z(s)\Big)\,ds\Big].
\end{split}
\end{equation}

And also
\begin{equation*}
\begin{split}
&E_t\int_t^Td\, \left\{\Big[X(s)-\kappa \int_0^{\delta}\!\!\!X(s\!-\!r)\mu_1(dr)\Big]Y(s)\right\}\\
=\,&E_t\Big[\Big(X(T)-\kappa\int_0^{\delta}\!\!\!X(T\!-\!r)\mu_1(dr)\Big)\,\xi\Big]-xY(t),
\end{split}
\end{equation*}
 In view of \eqref{1}, we obtain the conclusion.
\end{proof}

Then we introduce the Pontryagin maximum principle of the optimal
control Problem (C).
\begin{thm}
Suppose that $(H6)$ holds and U is convex. Let
$(\bar{X}(\cdot),\bar{u}(\cdot))$ be the optimal pair. Define
 $$\bar{b}_x(t):=
b_x(t,\bar{X}(t),\int_0^{\delta}\!\!\!\bar{X}(t\!-\!r)\mu_2(dr),\bar{u}(t)),$$
$$\bar{b}_y(t):=
b_y(t,\bar{X}(t),\int_0^{\delta}\!\!\!\bar{X}(t\!-\!r)\mu_2(dr),\bar{u}(t)),$$
$$\bar{b}_u(t):=b_u(t,\bar{X}(t),\int_0^{\delta}\!\!\!\bar{X}(t\!-\!r)\mu_2(dr),\bar{u}(t)),$$
$$\bar{l}_x(t):=
l_x(t,\bar{X}(t),\bar{u}(t)),~~~~~~~\bar{l}_u(t):=
l_u(t,\bar{X}(t),\bar{u}(t)).$$ where $b_y$ denotes the derivative
of $\,b$ on the third variable. Let $Y(\cdot)$ satisfy the following
equation:
\begin{equation*}
\left\{
\begin{split}
&-d\Big[Y(t)-\kappa
         E_t\int_0^{\delta}\!\!\!Y(t\!+\!r)\mu_1(dr)\Big]=\Big[\bar{b}_x(t)Y(t)
            +E_t\int_0^{\delta}\!\!\!\bar{b}_y(t\!+\!r)Y(t\!+\!r)\mu_2(dr)\\
&~~~~~~~~~~~~~~~~~~~~~~~~~~~~~~~~~~~+c(t,\bar{u}(t))Z(t)+\bar{l}_x(t)\Big]\,dt-Z(t)\,dW(t),~~~~~~t\in[0,T];\\
&Y(T)=M_x(\bar{X}(T));\\
&Y(t)=0,~~~~t\in(T,T+\delta].
\end{split}
\right.
\end{equation*}
Then we have almost surely
$$\big[\bar{b}_u(t)Y(t)+c_u(t,\bar{u}(t))\bar{X}(t)Z(t)+\bar{l}_u(t)\big](u-\bar{u}(t))\leq0,~~\forall u\in U.$$
\end{thm}

\begin{proof}
Since $U$ is convex, for all $u(\cdot)\in\sU$
$$u^{\eps}(\cdot)=\eps u(\cdot)+(1-\eps)\bar{u}(\cdot)\in\sU.$$ Let
$X^{\eps}$ and $\bar{X}$ be solutions of Eq. \eqref{NDSde}
corresponding to $u^{\eps}(\cdot)$ and $\bar{u}(\cdot)$,
respectively.

Define $$\Delta
X^{\eps}(t):=\frac{X^{\eps}(t)-\bar{X}(t)}{\eps},~~~~\Delta
u(t):=u(t)-\bar{u}(t),$$ then $\Delta X^{\eps}$ converges to $\chi$
in $\sS^2([-\delta,0],\bR)$, where $\chi$ satisfying
\begin{equation*}
\left\{
\begin{split}
&d\Big[\chi(t)-\kappa\int_0^{\delta}\!\!\!\chi(t\!-\!r)\mu_1(dr)\Big]
     =\Big[\bar{b}_x(t)\chi(t)+\bar{b}_y(t)\int_0^{\delta}\!\!\!\chi(t\!-\!r)\mu_2(dr)
           +\bar{b}_u(t)\Delta u(t)\Big]\,dt\\
&~~~~~
    +\!\!\Big[c(t,\bar{u}(t))\chi(t)\!+\!\theta(t)\!-\!\kappa\!\int_0^{\delta}\!\!\!
     \Big(c(t\!-\!r,\bar{u}(t\!-\!r))\chi(t\!-\!r)\!+\!\theta(t\!-\!r)\Big)\mu_2(dr)\Big]\,dW(t),~~t\in[0,T];\\
&\chi(t)=0,~~~~~~t\in[-\delta,0]
\end{split}
\right.
\end{equation*}
with $\theta(t):= c_u(t,\bar{u}(t))\bar{X}(t)\Delta u(t)$.

From the preceding lemma, we have
\begin{equation*}
\begin{split}
0\geq\,&\frac{J(u(\cdot))-J(\bar{u}(\cdot))}{\eps}\\
    =\,&E\left[\int_0^T\frac{l(t,X^{\eps}(t),u^{\eps}(t))-l(t,\bar{X}(t),\bar{u}(t))}{\eps}\,dt
          +\frac{M(T,X^{\eps}(T))-M(T,X^{\eps}(T))}{\eps}\right]\\
\rrow\,&E\left[\int_0^T\Bigl(\bar{l}_x(t)\chi(t)+\bar{l}_u(t)\Delta u(t)\Bigr)\,dt+\bar{M}_x\chi(T)\right]\\
    =\,&E\left[\int_0^T\Big(\bar{b}_u(t)\Delta u(t)Y(t)+\theta(t)Z(t)+\bar{l}_u(t)\Delta u(t)\Big)\,dt\right].
\end{split}
\end{equation*}

So we have almost surely
$$\big[\bar{b}_u(t)Y(t)+c_u(t,\bar{u}(t))\bar{X}(t)Z(t)+\bar{l}_u(t)\big](u-\bar{u}(t))\leq0,~\forall u\in U.$$
\end{proof}

\subsection{The optimal value}
In this subsection, we discuss the state-linear case, and construct
an explicit expression of its optimal value via the comparison
theorem.

Let $U$ be any nonempty subset of $\bR^m$. For all $t\in[0,T]$,
consider the following linear controlled system:
\begin{equation}\label{ForwardlDE}
\left\{
\begin{split}
&d\Big[X(s)-\kappa\int_0^{\delta}\!\!\!X(s\!-\!r)\mu_1(dr)\Big]=[a(s,u(s))X(s)
        +b(s,u(s))\int_0^{\delta}\!\!\!X(s\!-\!r)\mu_2(dr)]\,ds\\
&~~~~~~~~~~~~+\Big[c(s,u(s))X(s)-\kappa\int_0^{\delta}\!\!\!c(s\!-\!r,u(s\!-\!r))X(s\!-\!r)\mu_1(dr)\!\Big]\,dW(s),
          ~~s\in[t,T];\\
&X(t)=1;\\
&X(s)=0,~~~~~s\in[t-\delta,t),
\end{split}
\right.
\end{equation}
where $\kappa\in[0,1)$ is a constant,
$a(\cdot,u),b(\cdot,u)\in\sL^{\infty}_{\bF}(0,T;\bR^+)$ and
$c(\cdot,u)\in\sL^{\infty}_{\bF}(-\delta,T;\bR^d)$, $\forall u\in
U$.

The cost functional is also linear:
$$J(u(\cdot))(t)=E_t\left[\int_t^Tl(s,u(s))X(s)\,ds+MX(T)\right],$$
where $M\in\bL^2(\sF_T;\bR)$, and
$l(\cdot,u)\in\sL^2_{\bF}(0,T;\bR)$, $\forall u\in U$.

For all $t\in[0,T]$, the dynamic optimal control problem is the
following:

\textbf{Problem (C1)} Find a control process $\bar{u}(\cdot)\in\sU$
such that
$$J(\bar{u}(\cdot))(t)=\esssup_{u(\cdot)\in\sU}J(u(\cdot))(t).$$
$V(t):= J(\bar{u}(\cdot))(t)$ is the optimal value of Problem (C1).

In view of Lemma \ref{duality}, if $Y^u$ is the solution of the
following linear NBSFDE:
\begin{equation}\label{BlFE}
\left\{\!
\begin{split}
&-d\Big[Y^u(s)-\kappa E_s\int_0^{\delta}\!\!\!Y^u(s\!+\!r)\mu_1(dr)\Big]=\Big[a(s,u(s))Y^u(s)+c(s,u(s))Z^u(s)+l(s,u(s))\\
&~~~~~~~~~~~~~~~~~~~~~+E_s\int_0^{\delta}\!\!\!b(s\!+\!r,u(s\!+\!r))Y^u(s\!+\!r)\mu_2(dr)\Big]\,ds-Z(s)\,dW(s),~s\in[t,T];\\
&Y^u(T)=M;\\
&Y(s)=0,~~~~s\in(T,T+\delta],
\end{split}
\right.
\end{equation}
then $Y^u(t)=J(u(\cdot))(t)$.

For all $u\in\sU$, define
\begin{equation*}
\begin{split}
&f^u(s,Y(s),Y_s,Z(s)):= a(s,u(s))Y(s)+c(s,u(s))Z(s)+l(s,u(s))\\
&~~~~~~~~~~~~~~~~\,+E_s\Big[\int_0^{\delta}\!\!\!b(s\!+\!r,u(s\!+\!r))Y(s\!+\!r)\mu_2(dr)\Big].
\end{split}
\end{equation*} Then we have the following proposition which gives an
explicit expression of the optimal value of problem (C1).

\begin{prop}For all
$(y,z)\in\sH^2(0,T+\delta)$, set
$$\bar{f}(t,y(t),y_t,z(t))=\esssup_{u(\cdot)\in\sU}f^u(t,y(t),y_t,z(t)),$$
Then the NBSFDE
\begin{equation}\label{2}
\left\{
\begin{split}
&-d[Y(s)-\kappa
E_s\int_0^{\delta}\!\!\!Y(s\!+\!r)\mu_1(dr)]=\bar{f}(s,Y(s),Y_s,Z(s))\,ds-Z(s)\,dW(s),~~s\in[t,T];\\
&Y(T)=M;\\
&Y(s)=0,~~~~s\in(T,T+\delta]
\end{split}
\right.
\end{equation}
admits a unique adapted solution
$(\bar{Y},\bar{Z})\in\sH^2(0,T+\delta)$, and $\bar{Y}(t)=V(t)$.
\end{prop}

\begin{proof}
Applying Theorem \ref{ex_uni}, it is obvious that Eq.\eqref{2}
admits a unique pair of solution
$(\bar{Y},\bar{Z})\in\sH^2(0,T+\delta)$.

For all $u(\cdot)\in\sU$,
$$\bar{f}(t,\bar{Y}(t),\bar{Y}_t,\bar{Z}(t))
\geq f^u(t,\bar{Y}(t),\bar{Y}_t,\bar{Z}(t)).$$ In view of Theorem
\ref{compar}, for all $u(\cdot)\in\sU$, we have almost surely
$$\bar{Y}(t)\geq Y^u(t),~~~\forall t\in[0,T].$$

For all $\eps>0$,
choose an admissible control $u^{\eps}(\cdot)\in\sU$, such that
$$\bar{f}(t,\bar{Y}(t),\bar{Y}_t,\bar{Z}(t))\leq
f^{u^\eps}(t,\bar{Y}(t),\bar{Y}_t,\bar{Z}(t))+\eps.$$ Applying
Corollary \ref{depend_gener}, we have $$E\big[\sup_{0\leq t\leq
T}|Y^{u^{\eps}}(t)-\bar{Y}(t)|^2\big]\rrow 0,~~\textrm{as}
~\eps\rrow0$$

Therefore,
$$\bar{Y}(t)=\esssup_{u(\cdot)\in\sU}J(u(\cdot))(t)=V(t).$$
\end{proof}

\begin{rmk}
The controlled NSFDEs we consider here is of a special form in which
the diffusion has a similar form as the left side of the equation.
The limitation is due to the duality representation of Lemma
\ref{duality}. The optimal control of more general NSFDEs is still
open which will be discussed in the future.
\end{rmk}

\bibliographystyle{siam}

\end{document}